# A limit lifting theorem for fibrations between bicategories

Patrick Nicodemus, patrick.nicodemus5@gmail.com

## Abstract

A standard result from the theory of Grothendieck fibrations states that if $p : \mathbb{E} \to \mathbb{B}$ is a fibration, then $\mathbb{E}$ has limits of shape $\mathcal{J}$ if $\mathbb{B}$ has limits of shape $\mathcal{J}$, the fibers of $\mathbb{E}$ have limits of shape $\mathcal{J}$, and the reindexing functors preserve these limits.

Such a result is a basic tool when computing limits in a category known to be the total category of a fibration, for example, one that is defined using the Grothendieck construction.

This paper states and proves a generalization of this theorem for fibrations between bicategories.

## Notation

In a bicategory $\mathbb{B}$, we write both vertical and horizontal composition in diagrammatic order, so that if $f : x \to y$ and $g : y \to z$ are 1-cells, then $f \cdot g : x \to z$. We use the notation $(- \cdot -)$ to refer to horizontal composition of both 1-cells and 2-cells, and left and right whiskering. Vertical composition of 2-cells $\alpha : f \Rightarrow g, \beta : g \Rightarrow h$ is denoted $\alpha \mid \beta$.

For a 1-cell $f : x \to y$, we write $\ell_f$ for the left unitor $1_x \cdot f \cong f$ and $r_f$ for the right unitor $f \cdot 1_y$; note that the terms "left" and "right" are opposite to what is usually written, in order to be consistent with the diagrammatic order of composition. For $f : w \to x, g : x \to y, h : y \to z$ we write $a_{f,g,h}$ for the associator $(f \cdot g) \cdot h \Rightarrow f \cdot (g \cdot h)$.

If $p : \mathbb{E} \to \mathbb{B}$ is a lax functor or pseudofunctor between bicategories, then the lax functoriality constraint $p(f \cdot g) \Rightarrow p(f) \cdot p(g)$ is denoted $p^2_{f,g}$ and the lax unity constraint $1_{p(x)} \Rightarrow p(1_x)$ is denoted $p^0_x$. If $x, y$ are 0-cells in $\mathbb{B}$, $p_{x,y} : \mathbb{E}(x, y) \to \mathbb{B}(px, py)$ is the 1-functor component of $p$ between these hom-categories.

## Quasi-raris

If $F : C \to D$ is a functor between categories, and for every $d$ in $D$ there exists a universal arrow $(G(d), \varepsilon : F(G(d)) \to d)$ - i.e., a terminal object in the slice category $F \downarrow d$ - then the function $d \mapsto G(d)$ extends to a functor $G : D \to C$ which is right adjoint to $F$. In particular, if one can choose $G(d)$ so that $F(G(d)) = d$ and $\varepsilon$ is the identity, then there exists a functor $G : D \to C$ which is a **right adjoint right inverse** to $F$, or a **rari**. In this section we establish one generalization of this theorem to bicategories and pseudofunctors. The theory of quasi-adjunctions is complicated and is more general than is necessary for what follows, so we develop a specialized theory of **quasi-raris**.

Let $F : C \to D$ be a functor between categories. Given $d$ in $D$, $c$ is called a **rari-universal lift of $d$ along $F$** if $F(c) = d$ and $(F(c), \mathrm{id}_c)$ is a universal arrow $F \downarrow d$. We can omit $d$ and simply say "$c$ is rari-universal with respect to $F$" if $c$ is a rari-universal lift of $F(c)$ along $F$, and if $F$ is evident from the context then we will simply say $c$ is rari-universal.

In what follows, we will frequently make use of the fact that if $c$ is rari-universal, then $F_{c',c} : \mathbf{Hom}_C(c', c) \to \mathbf{Hom}_D(F(c'), F(c))$ is a bijection for all $c'$.

**Definition 0.1** (2-rari-universal lift for bicategories): Let $p : \mathbb{E} \to \mathbb{B}$ be a strict functor between bicategories. Let $b$ be an object in $\mathbb{B}$.

Then $e$ in $\mathbb{E}$ is a **2-rari-universal lift** of $b$ if $p(e) = b$ and, for each $e'$ in $\mathbb{E}$, the local hom functor $p_{e',e}$ has a right adjoint right inverse.

**Theorem 0.1**: Let $p : \mathbb{E} \to \mathbb{B}$ be a strict functor between bicategories. If for each $b$ in $\mathbb{B}$ there is a 2-rari-universal lift $e$ of $b$ along $p$, then there is a lax functor $s : \mathbb{B} \to \mathbb{E}$ which is a section of $p$.

*Proof*: For each $b$ in $\mathbb{B}$, choose a universal lift $e$ of $b$; this determines the behavior of $s$ on objects.

For $b, b'$ in $\mathbb{B}$, let $s_{b,b'}$ be defined as the assumed right adjoint to the projection functor $p_{s(b),s(b')}$.

Let us construct the lax unit constraint $1_{s(b)} \Rightarrow s(1_b)$. Because $s_{b,b}$ is right adjoint to the projection functor, it suffices to give a 2-cell $p\left(1_{s(b)}\right) \Rightarrow 1_b$; we take the identity 2-cell.

Let us construct the lax functoriality constraint $s(f) \cdot s(g) \Rightarrow s(f \cdot g)$, for $f : a \to b, g : b \to c$. Because $s_{a,c}$ is right adjoint to the projection functor, it suffices to give a map $p(s(f) \cdot s(g)) \Rightarrow f \cdot g$. We take the identity 2-cell.

To verify that the lax functoriality constraint is natural in $f$ and $g$, so let $f, f' : a \to b$ and $g, g' : b \to c$ with $\alpha : f \Rightarrow f', \beta : g \Rightarrow g'$, one uses that for any 1-cell $y$, $p_{s(a),s(b)}(s(f \cdot g), y)$ is bijective between sets of 2-cells.

The coherence conditions for associativity and left/right unitors are easy to prove, given that $p$ is strict. ∎

**Theorem 0.2**: Let $p : \mathbb{E} \to \mathbb{B}$ be a strict functor between bicategories.

Call the **identity condition** the following property: for each $b$ in $\mathbb{B}$ and each 2-rari-universal lift $e$ of $b$, $1_e$ is a rari-universal lift of $1_b$ along $p_{e,e}$.

Call the **composition condition** the condition that for each $x, y, z$ in $\mathbb{E}$ and $f : x \to y, g : y \to z$, if $f$ is rari-universal over $p(f)$ and $g$ is rari-universal over $p(g)$, then $f \cdot g$ is rari-universal over $p(f) \cdot p(g)$.

Suppose that $p$ satisfies the identity and composition conditions. Then if $e$ and $e'$ are 2-rari-universal lifts of $b, b'$ respectively, an equivalence between $b$ and $b'$ lifts to an equivalence between $e$ and $e'$. In particular, 2-rari-universal lifts are unique up to equivalence when they exist.

Furthermore, under the hypotheses of Theorem 0.1, the section $s$ thus constructed is a pseudofunctor, and is unique up to an an equivalence which lives over the identity equivalence of $1_\mathbb{B}$ with itself, that is, given two pseudofunctors $s, s'$ where $s(x), s'(x)$ are 2-lari-universal over $x$ and $s(f), s'(f)$ are rari-universal over $f$, $s, s'$ are equivalent by an equivalence $\tau$ such that $\tau \cdot p = 1_{1_\mathbb{B}}$.

*Proof*: Let $e, e'$ be 2-rari-universal lifts of $b, b'$ respectively, with an equivalence $f : b \to b', g : b' \to b, \eta : 1_b \cong f \cdot g, \varepsilon : g \cdot f \cong 1_{b'}$.

Let $\hat{f}, \hat{g}$ be rari-universal lifts of $f, g$ respectively, then $\hat{f} \cdot \hat{g}$ is rari-universal over $f \cdot g$ and $1_e$ is rari-universal over $1_b$. Thus the isomorphism $\eta$ lifts to a unique isomorphism $1_e \cong \hat{f} \cdot \hat{g}$. Similarly there is a unique lift of $\varepsilon$ to an isomorphism $\hat{g} \cdot \hat{f} \cong 1_{e'}$.

Given $f : x \to y, g : y \to z$ in $\mathbb{B}$, since $s(f) \cdot s(g)$ and $s(f \cdot g)$ are both rari-universal over $f \cdot g$, they are isomorphic by a unique vertical isomorphism, i.e., the lax functoriality constraint map constructed in the proof of Theorem 0.1 is an isomorphism. The same is true for the left and right unitor.

Let $s, s'$ be pseudofunctors satisfying the conditions of the theorem. We define a 1-cell $\tau_b : s(b) \to s'(b)$ for each $b$ as the rari-universal lift of $1_b$. The pseudonatural transformation 2-cell $\tau(f) : s(f) \cdot \tau_{b'} \Rightarrow \tau_b \cdot s(f')$ is the unique lift of the canonical isomorphism $f \cdot 1_{b'} \Rightarrow 1_b \cdot f$.

It is straightforward to verify that $\tau_f$ is natural in $f$, and that the lax unity and lax naturality conditions are satisfied; each argument is by appeal to the fact that $p_{e,e'}(f, f')$ is bijective when $f'$ is rari-universal.

If $\tau, \tau'$ are two such transformations constructed according to the specification above so that $\tau_b$ and $\tau'_b$ are both rari-universal then by their universal property $\tau_b \cong \tau'_b$ by a unique vertical isomorphism, and this extends to a unique vertical modification. Thus, $s$ is unique up to equivalence. ∎

If $F : \mathbb{C} \to \mathbb{D}$ satisfies the conditions of Theorem 0.1 and Theorem 0.2, we call the essentially unique section $s$ a **quasi-rari** for $F$.

It is outside the scope of this paper to fully relate these quasi-raris to quasi-adjunctions. However, we will establish a basic link to the existing literature, relating our notion to the presentation in sections 7.2 and 7.3 of [1]. This theorem is not used in the remainder of the paper.

Let $F : \mathbb{C} \to \mathbb{D}$ be a pseudofunctor between bicategories. Let $\mathbb{D}(F, 1)$ be the bicategory whose objects are triples $(c, d, f : F(c) \to d)$, whose 1-cells $(c, d, f) \to (c', d', f')$ are lax squares $(g : c \to c', h : d \to d', \alpha : F(g) \cdot f' \Rightarrow f \cdot h$, and whose 2-cells $(g, h, \alpha) \Rightarrow (g', h', \alpha')$ are pairs $(\sigma : g \Rightarrow g', \rho : h \Rightarrow h')$ making everything commute. Then $\mathbb{D}(F, 1)$ is equipped with a strict projection functor $\pi : \mathbb{D}(F, 1) \to \mathbb{D}$ and it is of interest to elaborate the conditions described in Theorem 0.1 and Theorem 0.2 in this setting. Actually, we will consider the 2-cell dualization of $\pi$, $\pi^{\text{co}} : \mathbb{D}(F, 1)^{\text{co}} \to \mathbb{D}^{\text{co}}$ to better fit the conventions of [1].

Suppose that each $d$ in $\mathbb{D}^{\text{co}}$ has a rari-universal lift along $\pi^{\text{co}}$, as in the assumption of Theorem 0.1, so that for every $d$ there is given a pair $(c, f : F(c) \to c)$ such that for every object $(c', d', f')$ in $\mathbb{D}$, the associated functor $\pi_{(c',d',f'),(c,d,f)}$ has a lari. Then every lax slice bicategory $F \downarrow d$ is equipped with an **inc-lax transformation** of the identity pseudofunctor to the constant pseudofunctor at $(c, f)$ (see [1], 7.2).

If a rari-universal lift $(c, f : F(c) \to d)$ of $d$ along $\pi^{\text{co}}$ additionally has the property specified in Theorem 0.2 that $(1_c, 1_d)$ is rari-universal, then $(c, f)$ is an **inc-lax terminal object** in $F \downarrow d$ (see [1], 7.2).

Last, if rari-universal 1-cells are closed under composition, then the change-of slice functors $F \downarrow u : F \downarrow d \to F \downarrow d'$ associated to a 1-cell $u : d \to d'$ **preserve inc-lax terminal objects** (see [1], 7.2).

This amounts to the following theorem:

> **Theorem 0.3**: Let $F : \mathbb{C} \to \mathbb{D}$ be a pseudofunctor between bicategories. Suppose that the functor $\pi^{\text{co}} : \mathbb{D}(F, 1)^{\text{co}} \to \mathbb{D}^{\text{co}}$ satisfies the hypotheses of Theorem 0.1 and Theorem 0.2. Then $F$ satisfies the hypotheses of the Bicategorical Quillen Theorem A from 7.3 of [1].

## Fibrations between bicategories

Fibrations of bicategories have been studied in [2] and [3]. We will follow [2] and recapitulate some of the definitions.

> **Definition 0.2** (Cartesian 1-cell): Let $\mathbb{E}, \mathbb{B}$ be bicategories, and $p : \mathbb{E} \to \mathbb{B}$ a lax functor. Let $f : x \to y$ be a 1-cell in $\mathbb{E}$. $f$ is **Cartesian** when the following properties are satisfied:
> 1. For each 0-cell $z$ in $\mathbb{E}$ and 1-cell $g : z \to y$, and for each 1-cell $h : pz \to px$ and 2-isomorphism $\alpha : pf \cdot h \Rightarrow pg$, there is a 1-cell $\hat{h} : z \to x$ and 2-isomorphisms $\hat{\alpha} : \hat{h} \cdot f \Rightarrow g, \hat{\beta} : p\hat{h} \Rightarrow h$ such that $\hat{\beta} \cdot pf \mid \alpha = p_{\hat{h},f}^2 \mid p\hat{\alpha}$. The triple $(\hat{h}, \hat{\alpha}, \hat{\beta})$ is called a **lift** of $(h, \alpha, g)$.
> 2. Let $z$ be a 0-cell in $\mathbb{E}$. Let $(g : z \to y, h : pz \to px, \alpha : h \cdot pf \cong pg)$ and $(g' : z \to x, h' : pz \to px, \alpha' : h' \cdot pf \cong pg')$ be two triples as above. Let $(\delta : h \Rightarrow h', \sigma : g \Rightarrow g')$ be a pair of 2-cells such that $\alpha \mid p\sigma = \delta \cdot pf \mid \alpha'$. Let $(\hat{h}, \hat{\alpha}, \hat{\beta})$ be a lift of $(h, \alpha, g)$ and $(\hat{h}', \hat{\alpha}', \hat{\beta}')$ be a lift of $(h', \alpha', g')$. Then there is a unique 2-cell $\hat{\delta} : \hat{h} \Rightarrow \hat{h}'$ such that $\beta \mid \delta = p\hat{\delta} \mid \beta'$ and $\hat{\delta} \cdot f \mid \hat{\alpha}' = \hat{\alpha}$.

> **Definition 0.3** (Cartesian 2-cell): Let $\mathbb{E}, \mathbb{B}$ be bicategories, and $p : \mathbb{E} \to \mathbb{B}$ a lax functor. Let $f, g : x \to y$ be 1-cells in $\mathbb{E}$ and $\sigma : f \Rightarrow g$ a 2-cell. Then $\sigma$ is **Cartesian** as a 2-cell if it is Cartesian in the usual sense with respect to $p_{x,y} : \mathbb{E}(x, y) \to \mathbb{B}(px, py)$.

> **Definition 0.4** (Locally fibered pseudofunctor): Let $\mathbb{E}, \mathbb{B}$ be bicategories, and $p : \mathbb{E} \to \mathbb{B}$ a lax functor. $p$ is **locally fibred** if for all 0-cells $x, y : \mathbb{E}$, $p_{x,y}$ is a Grothendieck fibration.

> **Definition 0.5** (Fibration of bicategories): Let $\mathbb{E}, \mathbb{B}$ be bicategories, and $p : \mathbb{E} \to \mathbb{B}$ a pseudofunctor. $p$ is a **fibration** when:
> 1. $p$ is locally fibred
> 2. for each 0-cell $y$ in $E$ and 1-cell $f : x \to py$, there exists a Cartesian 1-cell $h : z \to y$ with $ph = z$
> 3. the horizontal composite of two Cartesian 2-cells is again Cartesian

The following theorem strengthens condition 1 of Definition 0.2 by imposing the constraint that $\hat{\beta}$ must be the identity 2-cell.

**Theorem 0.4** (Strictification of the lifting condition): Let $p : \mathbb{E} \to \mathbb{B}$ be a locally fibred pseudofunctor. Then a 1-cell $f : x \to y$ in $\mathbb{E}$ is Cartesian if and only if has the following lifting property:

1. for every 0-cell $z$ in $\mathbb{E}$ and 1-cell $g : z \to y$, and for every 1-cell $h : pz \to px$ and 2-isomorphism, $\alpha : pf \cdot h \Rightarrow pg$ there is a 1-cell $\hat{h} : z \to x$ and a 2-isomorphism $\hat{\alpha} : \hat{h} \cdot f \Rightarrow g$ such that $p\hat{h} = h$ and $p^2_{\hat{h},f} \mid p\hat{\alpha} = \alpha$.
2. Let $z$ be a 0-cell in $\mathbb{E}$. Let $(g : z \to y, h : pz \to px, \alpha : h \cdot pf \cong pg)$ and $(g' : z \to x, h' : pz \to px, \alpha' : h' \cdot pf \cong pg')$ be two triples as above. Let $(\delta : h \Rightarrow h', \sigma : g \Rightarrow g')$ be a pair of 2-cells such that $\alpha \mid p\sigma = \delta \cdot pf \mid \alpha'$. Let $(\hat{h}, \hat{\alpha})$ be a strict lift of $(h, \alpha, g)$ (i.e, $\hat{\beta}$ is the identity 2-cell) and and $(\hat{h}', \hat{\alpha}')$ be a strict lift of $(h', \alpha', g')$. Then there is a unique 2-cell $\hat{\delta} : \hat{h} \Rightarrow \hat{h}'$ such that $\delta = p\hat{\delta}$ and $\hat{\delta} \cdot f \mid \hat{\alpha}' = \hat{\alpha}$.

*Proof*: Note that point (1.) is a **strengthening** of the Cartesian lifting condition for 1-cells, while (2.) is a **weakening** of the Cartesian lifting condition for 2-cells.

The forward implication (if $f$ is Cartesian, then it has the strict lifting property) is Proposition 3.2.1 of [2].

The reverse implication is by a similar proof, and we omit the details. ∎

**Theorem 0.5** (Dropping the invertibility of $\alpha$): Let $p : \mathbb{E} \to \mathbb{B}$ be a locally fibred pseudofunctor, and let $f$ be a Cartesian 1-cell. Then $f$ has the analogous lifting property given by dropping the assumption that $\alpha$ be invertible, in Definition 0.2., and weakening the conclusion "$\hat{\alpha}$ is a 2-isomorphism" to "$\hat{\alpha}$ is Cartesian". More explicitly, a Cartesian 1-cell has the following lifting properties:

1. for every 0-cell $z$ in $\mathbb{E}$ and 1-cell $g : z \to y$, and for every 1-cell $h : pz \to px$ and not-necessarily invertible 2-cell $\alpha : pf \cdot h \Rightarrow pg$ there is a 1-cell $\hat{h} : z \to x$ and a Cartesian 2-cell $\hat{\alpha} : \hat{h} \cdot f \Rightarrow g$ such that $p\hat{h} = h$ and $p^2_{\hat{h},f} \mid p\hat{\alpha} = \alpha$.
2. for every 0-cell $z$ in $\mathbb{E}$, and any two lifting problems $(g_1 : z \to y, h_1 : pz \to px, \alpha_1 : h_1 \cdot pf \Rightarrow pg_1)$ and $(g_2 : z \to y, h_2 : pz \to px, \alpha_2 : h_2 \cdot pf \Rightarrow pg_2)$; for any "morphism of lifting problems" $(\delta : h_1 \Rightarrow h_2, \sigma : g_1 \Rightarrow g_2)$ such that $\delta \cdot pf \mid \alpha_2 = \alpha_1 \mid p\sigma$ and any two lifts $(\hat{h}_1 : z \to x, \hat{\alpha}_1 : \hat{h}_1 \cdot f \Rightarrow g_1, \hat{\beta}_1 : p\hat{h}_1 \cong h_1)$ and $(\hat{h}_2 : z \to x, \hat{\alpha}_2 : \hat{h}_2 \cdot f \Rightarrow g_2, \hat{\beta}_2 : p\hat{h}2 \cong h_2)$, if $\hat{\alpha}_2$ is Cartesian, there is a unique 2-cell $\hat{\delta} : \hat{h}_1 \Rightarrow \hat{h}_2$ such that $p(\hat{\delta}) \mid \hat{\beta}_2 = \hat{\beta}_1 \mid \delta$ and $\hat{\delta} \cdot f \mid \hat{\alpha}_2 = \hat{\alpha}_1 \mid \sigma$.

*Proof*:
1. Let $\rho$ be a Cartesian 2-cell over $\alpha$ with codomain $g$, and write $k$ for its domain. Then by the lifting property of a Cartesian 1-cell, there is a 1-cell $\hat{h}$ and a 2-isomorphism $\gamma : \hat{h} \cdot f \cong k$ such that $p\hat{h} = h$ and $p\gamma = p^2_{\hat{h},f}$. Then we define $\hat{\alpha} = \gamma \mid \rho$. $\rho$ was Cartesian by assumption, $\gamma$ is an isomorphism, and the composition of Cartesian 2-cells is Cartesian, so we are done.
2. Apply the 2-lifting property of a Cartesian 2-cell to $\hat{\alpha}_2$; the map $p(\hat{\alpha}_1 \mid \sigma)$ factors through $p\hat{\alpha}_2$ according to
$$p^2_{\hat{h}_1,f}{}^{-1} \mid \hat{\beta}_1 \mid \delta \cdot f \mid p^2_{\hat{h}_2,f} \mid p\hat{\alpha}_2 = p(\hat{\alpha}_1 \mid \sigma)$$
Therefore there is a unique 2-cell $\kappa : \hat{h}_1 \cdot f \Rightarrow \hat{h}_2 \cdot f$ such that

$$p\kappa = {p^2_{\hat{h}_1, f}}^{-1} \mid \hat{\beta}_1 \mid \delta \cdot f \mid p^2_{\hat{h}_2, f}$$

and $\kappa \mid \hat{\alpha}_2 = \hat{\alpha}_1 \mid \sigma$. We can apply the known 2-lifting property of the Cartesian 1-cell $f$ to the lifting problems $\left(h_1, p^2_{\hat{h}_1, f}, \hat{h}_1 \cdot f\right)$ and $\left(h_2, p^2_{\hat{h}_2, f}, \hat{h}_2 \cdot f\right)$ with lifts $\left(\hat{h}_1, 1_{\hat{h}_1 \cdot f}, \hat{\beta}_1\right)$ and $\left(\hat{h}_2, 1_{\hat{h}_2 \cdot f}, \hat{\beta}_2\right)$ respectively, with $(\delta, \kappa)$ the 2-cells connecting the lifting problems. One easily sees that $p^2_{\hat{h}_1, f} \mid p\kappa = \delta \cdot pf \mid p^2_{\hat{h}_2, f}$ so that we can apply lifting condition. Thus there is a unique 2-cell $\hat{\delta}$ such that $\hat{\beta}_1 \mid \delta = p\hat{\delta} \mid \hat{\beta}_2$ (which is the condition we want) and $\hat{\delta} \cdot f = \kappa$, so $\hat{\delta} \cdot f \mid \hat{\alpha}_2 = \hat{\alpha}_1 \mid \sigma$. If $\hat{\delta}'$ is any morphism satisfying these same properties, then by 2-Cartesianity of $\hat{\alpha}_2$ we must have $\hat{\delta}' \cdot f = \kappa$, so $\hat{\delta} = \hat{\delta}'$.

∎

**Theorem 0.6** (Fibrations can be made strict): Let $p : \mathbb{E} \to \mathbb{B}$ be a locally fibered pseudofunctor. Then in the slice 1-category of pseudofunctors $p : \mathbb{X} \to \mathbb{B}$ and pseudofunctors $s : (\mathbb{X}, p) \to (\mathbb{Y}, q)$ such that $q = s \cdot p$ on the nose, $p$ is isomorphic to a **strict** locally fibered pseudofunctor $p'$.

*Proof*: This is Proposition 3.2.2 of [2]. ∎

Henceforth, all locally fibred pseudofunctors we introduce are assumed to be strict, and when we appeal to the lifting property of a Cartesian 1-cell in a fibration, we make use of the stronger property given in Theorem 0.4.

**Definition 0.6** (Cloven fibration): A fibration $p : \mathbb{E} \to \mathbb{B}$ between bicategories is **cloven** when:
- for each $e$ in $E$ and $f : b \to pe$, there is given a choice of Cartesian lift of $f$ with codomain $e$
- for each Cartesian 1-cell $f$ in $\mathbb{E}$, and each triple $(h, \alpha, g)$ as in Definition 0.2, there is given a choice of lift $\left(\hat{h}, \hat{\alpha}, \hat{\beta}\right)$.
- for all $x, y \in \mathbb{E}$, the hom-functor $p_{x,y} : \mathbb{E}(x, y) \to \mathbb{B}(x, y)$ is a cloven fibration

Whenever we assume a fibration is cloven, we will further assume that the map $(h, \alpha, g) \mapsto \left(\hat{h}, \hat{\alpha}, \hat{\beta}\right)$ is chosen such that $\hat{\beta}$ is always the identity 2-cell.

Let us isolate a more convenient form of the uniqueness condition that reflects the way we will use it.

**Theorem 0.7**: Let $p : \mathbb{E} \to \mathbb{B}$ be a pseudofunctor. Let $z, x \in \mathbb{E}$, let $h, h' : z \to x$, and let $\hat{\delta}_1, \hat{\delta}_2 : h \Rightarrow h'$. Let $f : x \to y$ be any Cartesian 1-cell, and let $\alpha : h' \cdot f \Rightarrow g$ be any Cartesian 2-cell. Then to prove $\hat{\delta}_1 = \hat{\delta}_2$, it suffices to show that $p\bigl(\hat{\delta}_1\bigr) = p\bigl(\hat{\delta}_2\bigr)$ and that $\hat{\delta}_1 \cdot f \mid \alpha = \hat{\alpha}_2 \cdot f \mid \alpha$.

*Proof*: Any Cartesian 2-cell is also monic, so $\hat{\alpha}_1 \cdot f \mid \alpha = \hat{\alpha}_2 \cdot f \mid \alpha$ implies $\hat{\alpha}_1 \cdot f = \hat{\alpha}_2 \cdot f$. One then can apply the unique factoring property of the Cartesian 1-cell $f$, taking $\sigma = \hat{\alpha}_1 \cdot f \mid \alpha = \hat{\alpha}_2 \cdot f \mid \alpha$. ∎

# The free fibration

Let $p : E \to B$ be a 1-functor between 1-categories. There is an associated "free fibration" over $B$ generated by $p$, defined as follows:
- the objects are triples $(b, e, f : b \to pe)$
- a morphism $(b, e, f) \to (b', e', f')$ is a pair of morphisms $(s : b \to b', t : e \to e')$ such that $s \cdot f' = f \cdot pt$.

This construction can be extended to a 2-functor on $\mathbf{Cat}/B$; it is a Kock–Zöberlein 2-monad on $\mathbf{Cat}/B$ whose algebras are precisely fibrations equipped with a cleavage (see [1], section 9.2).

Something similar is true for fibrations between bicategories. ([2], section 4.2) and [3] give a construction of the free fibration on a pseudofunctor between bicategories; we recall this here.

**Definition 0.7** (Oplax comma category associated to a pseudofunctor): Let $p : \mathbb{E} \to \mathbb{B}$ be a pseudofunctor between bicategories. The bicategory $\mathbb{B}/p$ is defined as follows:
- its objects are ordered triples $(b, e, f : b \to pe)$, where $b \in \mathbb{B}_0$ and $e \in \mathbb{E}_0$
- A 1-cell $(b, e, f) \to (b', e', f')$ is defined to be a triple $(s, t, \alpha)$, where $s : b \to b'$, $t : e \to e'$, and $\alpha : s \cdot f' \Rightarrow f \cdot pt$.
- A 2-cell $(s, t, \alpha) \Rightarrow (s', t', \alpha')$ consists of a pair $(\beta : s \Rightarrow s', \gamma : t \Rightarrow t')$ such that $\beta \cdot f' \mid \alpha' = \alpha \mid f \cdot p\gamma$
- composition and identity are evident
- the associator and unitors, and the coherence conditions, are inherited from $\mathbb{B}$ and $\mathbb{E}$

$\mathbb{B}/p$ is equipped with a strict projection functor $d_0 : \mathbb{B}/p \to \mathbb{B}$ mapping $(b, e, f) \mapsto b$ and $(s, t, \alpha) \mapsto s$.

**Theorem 0.8**: For any pseudofunctor $p : \mathbb{E} \to \mathbb{B}$, $d_0$ is a fibration.

*Proof*: This is proposition 4.2.4 of [2] and Theorem 4.24 of [3]. ∎

# Cartesian lifts as universal arrows

Let $p : \mathbb{E} \to \mathbb{B}$ be a locally fibred pseudofunctor, and let $\mathbb{E}^I$ denote the bicategory of arrows and oplax squares in $\mathbb{E}$. Then $p$ induces a pseudofunctor $p^L : \mathbb{E}^I \to \mathbb{B}/p$ sending $f : e \to e'$ to $(p(b), e', f)$.

**Theorem 0.9**: Let $(b, e, f : b \to p(e))$ be an object in $\mathbb{B}/p$ and let $\hat{f}$ be a Cartesian lift of $f$. Then $\hat{f}$ is a 2-rari-universal lift of $(b, e, f)$ along $p^L$.

*Proof*: Let $w_0 : e_0 \to e_0'$ be an object in $\mathbb{E}^I$ and let $(s : p(e_0) \to p(e_1), t : e_0' \to e_1', \alpha) : p^{L(w_0)} \to p^{L(w_1)}$ in $\mathbb{B}/p$ from $p^L(w_0) \to p^L(w_1)$.

To show that $\left(p^L \downarrow p^{L(w_1)}\right)\left((w_0, (s, t, \alpha)), \left(w_1, 1_{p^{L(w_1)}}\right)\right)$ has a terminal object, it suffices to show that there is a 1-cell $\overline{\tau}(s, t, \alpha)$ in $\mathbb{E}^I$ lying strictly over $(s, t, \alpha)$, such that the identity 2-cell $p^L(\overline{\tau}(s, t, \alpha)) \Rightarrow (s, t, \alpha)$ is a universal arrow from $p^L_{w_0, w_1}$ to $(s, t, \alpha)$. This is a small reformulation of Theorem 0.5 and we omit the details. ∎

**Lemma 0.9.1**: $p^L$ always satisfies the identity condition of Theorem 0.2, and it satisfies the composition condition when the horizontal composite of Cartesian 2-cells in $\mathbb{E}$ is again Cartesian.

*Proof*: Straightforward. ∎

Therefore, we conclude:

**Theorem 0.10**: Let $p$ be a strict fibration. Then $p^L$ has an essentially unique quasi-rari $s^L$.

# The exponential fibration

Let $\mathbb{A}, \mathbb{E}$ be bicategories; then there is a bicategory $\mathrm{Bicat}^{\mathrm{op}}(\mathbb{A}, \mathbb{E})$ whose 0-cells are lax functors $\mathbb{A} \to \mathbb{E}$, whose 1-cells are oplax natural transformations, and whose 2-cells are modifications. (see [1], Theorem 4.4.11)

For a pseudofunctor $p : \mathbb{E} \to \mathbb{B}$, postwhiskering with $p$ induces a pseudofunctor $\mathrm{Bicat}^{\mathrm{op}}(\mathbb{A}, \mathbb{E}) \to \mathrm{Bicat}^{\mathrm{op}}(\mathbb{A}, \mathbb{B})$, which is strict if $p$ is (see [1], sections 11.1 and 11.3).

The main theorem of this section is as follows:

**Theorem 0.11** (Exponential fibration): Let $p : \mathbb{E} \to \mathbb{B}$ be a cloven fibration between bicategories. Then for any bicategory $\mathbb{A}$, the postwhiskering pseudofunctor $- \cdot p :$ $\mathrm{Bicat}^{\mathrm{op}}(\mathbb{A}, \mathbb{E}) \to \mathrm{Bicat}^{\mathrm{op}}(\mathbb{A}, \mathbb{B})$ is a fibration.

We will prove this in a series of lemmas. We assume that $p$ is strict, and make use of this fact without comment.

**Theorem 0.12**: Let $\mathbb{A}$ be a bicategory, and $p : \mathbb{E} \to \mathbb{B}$ a strict locally fibered pseudofunctor. Suppose that the horizontal composite of Cartesian 2-cells in $p$ is Cartesian. Then $(- \cdot p) :$ $\mathrm{Bicat}^{\mathrm{op}}(\mathbb{A}, \mathbb{E}) \to \mathrm{Bicat}^{\mathrm{op}}(\mathbb{A}, \mathbb{B})$ is locally fibered.

*Proof*: Let $F, G : \mathbb{A} \to \mathbb{E}$ be lax functors. Let $\tau : F \Rightarrow G$ be an oplax natural transformation. Let $\sigma : F \cdot p \Rightarrow G \cdot p$ be an oplax natural transformation. Let $\gamma : \sigma \Rightarrow \tau \cdot p$ be a modification. We want to prove that there is a lift $\overline{\sigma}$ of $\sigma$ along $p$ and a Cartesian lift $\overline{\gamma}$ of $\gamma$ along $p$.

For $a$ in $\mathbb{A}$, define $\overline{\gamma}(a)$ to be any choice of Cartesian 2-cell with codomain $\tau(a)$ such that $p(\overline{\gamma}(a))$; define $\overline{\sigma}(a)$ to be the domain of $\overline{\gamma}(a)$.

For $f : a \to a'$ in $\mathbb{A}$, define $\overline{\sigma}(f)$ as follows: the horizontal composite $\overline{\gamma}(a) \cdot G(f)$ is Cartesian because $\overline{\gamma}(a)$ is and $1_{G(f)}$ is. Thus to construct a 2-cell $\overline{\sigma}(f) : F(f) \cdot \overline{\tau}(a') \Rightarrow \overline{\tau}(a) \cdot G(f)$ it suffices to observe that $\sigma : p(F(f) \cdot \overline{\tau}(a')) \Rightarrow p(\overline{\tau}(a) \cdot G(f))$ makes the appropriate triangle commute in $\mathbb{B}(p(F(a)), p(G(a')))$ (because $\gamma$ is a modification).

It is necessary to verify that $\overline{\sigma}(f)$ is natural in $f$. The argument is by appeal to the fact that $\overline{\gamma}(a) \cdot G(f)$ is monic for any $f$, appeal to the definition of $\overline{\sigma}(f)$, and the naturality condition of $\tau$.

We also sketch the proof of the oplax unity and composition constraints. The proof technique is similar to the naturality of $\overline{\sigma}(f)$. As before we should use the fact that $\overline{\gamma}(a) \cdot G(1_f)$ is monic for any $f$ because it is Cartesian, and then reduce the unity / composition constraint for $\overline{\sigma}$ to the unity / composition constraint for $\tau$.

Thus, $\overline{\sigma}$ is an oplax natural transformation. It now follows immediately from the definition of $\overline{\gamma}$ that it is a modification.

It remains to prove that $\overline{\gamma}$ is Cartesian. Therefore, let $\rho$ be some oplax natural transformation $F \Rightarrow G$ and let $\theta : \rho \Rightarrow \tau$ be a modification; suppose that $\kappa : \rho \cdot p \Rightarrow \sigma$ is a modification such that $\kappa \mid \gamma = \theta \cdot p$. Then for each $a$ there is a unique lift $\overline{\kappa}(a)$ of $\kappa(a)$ along $p_{F(a),G(a)}$ such that $\overline{\kappa}(a) \mid \overline{\tau}(a) = \theta(a)$. To see that $\kappa$ is indeed a modification, we use the same proof technique as above, appealing to the fact that $\overline{\gamma}(a) \cdot G(f)$ is monic because it is Cartesian, and then appealing to the definition of $\overline{\gamma}$ and the fact that $\theta$ is a modification. ∎

> **Theorem 0.13**: Let $\mathbb{A}$ be a bicategory, and $p : \mathbb{E} \to \mathbb{B}$ a pseudofunctor.
> - There is a bijective correspondence between lax (respectively colax, pseudo) functors $\mathbb{A} \to \mathbb{B}/p$ and triples $(F, G, \tau)$, where $F : \mathbb{A} \to \mathbb{B}$ and $G : \mathbb{A} \to \mathbb{E}$ are lax (respectively colax, pseudo) functors and $\tau : F \Rightarrow G \cdot p$ is an oplax natural transformation.
> - There is a bijective correspondence between lax (respectively colax, pseudo) functors $\mathbb{A} \to \mathbb{E}^I$ and triples $(\overline{F}, G, \overline{\tau})$ where $\overline{F}$ and $G$ are lax (respectively colax, pseudo) functors $\mathbb{A} \to \mathbb{E}$ and $\overline{\tau}$ is an oplax natural transformation $\overline{F} \Rightarrow G$.
> - If $(F, G, \tau) : \mathbb{A} \to p/\mathbb{B}$ is given, then a lift (on the nose) of $(F, G, \tau)$ along $p^L$ precisely corresponds to a pair $(\overline{F}, \overline{\tau})$ where $\overline{F}$ is a lift of $F$ along $p$ and $\overline{\tau}$ is a lift of $\tau$ along $p$.

All of this is seen easily by inspecting the definitions.

> **Theorem 0.14** (Canonical lift associated to an oplax natural transformation in a fibration):
> Let $p : \mathbb{E} \to \mathbb{B}$ be a strict fibration, and $\mathbb{A}$ a bicategory. Let $F : \mathbb{A} \to \mathbb{B}$ and $G : \mathbb{A} \to \mathbb{E}$ be lax (respectively colax, pseudo) functors, and $\tau : F \Rightarrow G \cdot p$ an oplax natural transformation. Then there is a lax (respectively colax, pseudo) functor $\overline{F} : \mathbb{A} \to \mathbb{E}$ which is a strict lift of $F$ along $p$, and an oplax natural transformation $\overline{\tau} : \overline{F} \Rightarrow G$ such that the 1-cells and 2-cells of $\overline{\tau}$ are Cartesian.

*Proof*: It suffices to give a lift of $(F, G, \tau)$ along $p^L$; therefore we define the pair $(\overline{F}, \overline{\tau})$ as $(F, G, \tau) \cdot s^L$. It is immediate by construction that $\overline{F}$ and $\overline{\tau}$ are well defined and strict lifts of $F, \tau$ respectively, and that the 1- and 2-cells of $\overline{\tau}$ are Cartesian. ∎

> **Theorem 0.15**: Let $p : \mathbb{E} \to \mathbb{B}$ be a strict fibration, and $\mathbb{A}$ a bicategory. Then the postcomposition functor $- \cdot p : \text{Bicat}^{\text{op}}(\mathbb{A}, \mathbb{E}) \to \text{Bicat}^{\text{op}}(\mathbb{A}, \mathbb{B})$ has the Cartesian lifting property; a 1-cell in $\text{Bicat}^{\text{op}}(\mathbb{A}, \mathbb{E})$ is Cartesian if its 1-cells and 2-cells as an oplax natural transformation are Cartesian.
>
> The same is true for the associated bicategories of colax functors and oplax natural transformations, and pseudofunctors and oplax natural transformations.

*Proof*: Because we have established that $- \cdot p$ is locally fibered, we use the stricter characterization of Cartesian 1-cells in Theorem 0.4.

We only treat the case of lax functors, as the other cases are similar. For oplax functors, everything is the same except the oplax unity and composition constraints for a natural transformation. We regard a pseudofunctor as equipped with both a lax and oplax functor structure which are mutually inverse, so once the lax and oplax cases are proved it is only a matter of checking that the two lifted lax and oplax structures of a pseudofunctor are inverse.

Introduce lax functors $F : \mathbb{A} \to \mathbb{B}, G : \mathbb{A} \to \mathbb{E}$, and $\tau : F \Rightarrow G \cdot p$ an oplax natural transformation. We use the definitions of $\overline{F}$ and $\overline{\tau}$ given above.

It must be shown that $\overline{\tau}$ is a Cartesian 1-cell; we prove the two components of the definition separately.

Thus, introduce $H : \mathbb{A} \to \mathbb{E}$ a lax functor, and $\sigma : H \Rightarrow G$ an oplax natural transformation; introduce $\kappa : H \cdot p \Rightarrow F$ an oplax natural transformation, and $\alpha : \kappa \cdot \tau \cong \sigma \cdot p$ an isomorphism (i.e., an invertible modification); we will prove that there is a strict lift $\overline{\kappa} : H \Rightarrow \overline{F}$ of $\kappa$ along $p$ and a strict lift $\overline{\alpha} : \overline{\kappa} \cdot \overline{\tau} \cong \sigma$ of $\alpha$ along $p$.

For $a$ in $\mathbb{A}$, we define $\overline{\kappa}(a)$ and $\overline{\alpha}(a)$ in the evident way by appeal to Theorem 0.4.

For $f : a \to a'$ in $\mathbb{A}$, we define $\overline{\kappa}(f)$ as the unique 2-cell over $\kappa(f)$ such that the already given definition for $\overline{\alpha}(a)$ will constitute a modification $\overline{k} \cdot \overline{\tau} \cong \sigma$. That a unique such morphism exists follows from the fact that $\overline{\tau}(a')$ is Cartesian, and the pasting of the 2-cells $\overline{F}(f)$ with $\alpha(a)$ is Cartesian.

Let us prove that $\overline{\kappa}(f)$ is actually an oplax natural transformation, by verifying that $\overline{\kappa}$ is natural in $f$, and that the oplax unity and composition constraints hold. Each is an application of Theorem 0.7 and the proof techniques are very similar.

Let $f, g : a \to a'$ in $\mathbb{A}$, and $\rho : f \Rightarrow g$. To establish the naturality of $\overline{\kappa}$ for $\rho$, apply Theorem 0.7 with Cartesian 1-cell $\overline{\tau}(a')$ and Cartesian 2-cell $\overline{\alpha}(a) \cdot \overline{\tau}(g)$; then the naturality of $\overline{\kappa}$ at $\rho$ follows from the naturality of $\overline{\tau}$ and $\sigma$ at $\rho$.

For $a$ in $\mathbb{A}$, the oplax unity condition at $a$ is established by applying Theorem 0.7 with Cartesian 1-cell $\overline{\tau}(a)$ and Cartesian 2-cell $\overline{\alpha}(a) \cdot \overline{\tau}(1_a)$, appealing to the oplax unity constraint for $\overline{F}$ and $G$.

For $f : a \to b, g : b \to c$ in $\mathbb{A}$, the oplax composition condition at $(f, g)$ is established by applying Theorem 0.7 with Cartesian 1-cell $\overline{\tau}(c)$ and Cartesian 2-cell $\overline{\alpha}(a) \cdot G(f \cdot g)$.

For part 2 of the definition of Cartesian, let $g, g' : H \Rightarrow G$ and $h, h' : H \Rightarrow \overline{F}$ be oplax natural transformations. Let $\sigma : g \Rightarrow g', \delta : h \cdot p \Rightarrow h' \cdot p, \alpha : h \cdot \tau \cong g, \alpha' : h' \cdot \tau \cong g'$, and $\delta : h \cdot p \Rightarrow h' \cdot p$ be modifications, such that $\delta \cdot \tau \mid \alpha' \cdot p = \alpha \cdot p$. We will show there is a unique modification $\hat{\delta} : h \Rightarrow h'$ such that $p(\hat{\delta}) = \delta$ and $\delta \cdot \overline{\tau} \mid \alpha' = \alpha$.

The modification is constructed pointwise in the obvious way, for each $a \in \mathbb{A}$, appealing to (2.) of , because each $\overline{\tau}(a)$ is Cartesian. It is evidently unique, so it suffices to check that it satisfies the axiom for a modification. Let $f : a \to a'$ in $\mathbb{A}$. We apply Theorem 0.4 with Cartesian 1-cell $\overline{\tau}(a)$ and Cartesian 2-cell $h'(a) \cdot \overline{\tau}(f) \mid \overline{\alpha}(a)$, then use the hypothesis that $\sigma : g \Rightarrow g'$ is a modification. ∎

**Theorem 0.16**: Let $p : \mathbb{E} \to \mathbb{B}$ be a strict fibration, and $\mathbb{A}$ a bicategory. A 1-cell $\tau : F \Rightarrow G$ in $\text{Bicat}^{\text{op}}$ ($\mathbb{A}, \mathbb{E}$) is Cartesian iff it has Cartesian 1-cells and Cartesian 2-cells as an oplax natural transformation.

*Proof*: We have proved one direction of the implication in Theorem 0.15. For the other, let $\tau : F \Rightarrow G$ be a Cartesian oplax natural transformation in $\text{Bicat}^{\text{op}}$ ($\mathbb{A}, \mathbb{E}$). Let $\overline{\tau}$ be the Cartesian oplax natural transformation of $\tau \cdot p$ with codomain $G$ we constructed in Theorem 0.15; then $\tau$ factors through $\overline{\tau}$ by some $\sigma : F \Rightarrow \overline{F}$ up to an invertible modification. Evidently $\sigma$ is a pseudonatural equivalence, because we can construct $\sigma^{-1}$ going in the opposite direction; thus, the 1-cells of $\sigma$ are equivalences and the 2-cells are isomorphisms.

By 3.1.8 and 3.1.9 of [2], it follows that the 1-cells of $\tau$ are Cartesian.

Similarly, because the pasting of Cartesian 2-cells is Cartesian, $\sigma(f)$ is an isomorphism (thus Cartesian) and $\overline{\tau}(f)$ is an isomorphism, thus the pasting is an isomorphism. ∎

**Theorem 0.17**: Let $p : \mathbb{E} \to \mathbb{B}$ be a strict fibration, and $\mathbb{A}$ a bicategory. The horizontal composite of Cartesian 1-cells in $\text{Bicat}^{\text{op}}$ ($\mathbb{A}, \mathbb{E}$) is again Cartesian.

*Proof*: We use Theorem 0.16. Then apply [2] 3.1.9 which states that the composition of Cartesian 1-cells is Cartesian. Because $p$ is a fibration, the pasting of Cartesian 2-cells is again Cartesian. ∎

This concludes the proof of Theorem 0.11.

**Theorem 0.18**: For bicategories $\mathbb{A}, \mathbb{B}$, let $\text{Bicat}^{\text{ps}}$ ($\mathbb{A}, \mathbb{B}$) denote the bicategory of lax functors (respectively, oplax and pseudo) and pseudonatural transformations.

Let $p : \mathbb{E} \to \mathbb{B}$ be a cloven fibration between bicategories.

Then for any bicategory $\mathbb{A}$, the postwhiskering pseudofunctor $- \cdot p : \text{Bicat}^{\text{ps}}$ ($\mathbb{A}, \mathbb{E}$) $\to \text{Bicat}^{\text{ps}}$ ($\mathbb{A}, \mathbb{B}$) is a fibration.

*Proof*: In a 1-Grothendieck fibration, a Cartesian lift of an isomorphism is an isomorphism. It follows that the functor $s^L$ sends 1-cells whose underlying square is invertible to 1-cells whose underlying square is invertible, and that a Cartesian lift of a pseudonatural transformation between lax functors is again a pseudonatural transformation.

The other thing that needs to be proven is that the factoring of a pseudonatural transformation through a Cartesian pseudonatural transformation is pseudo. Inspecting the definition of $\overline{\kappa}(f)$, it suffices to observe that the functor $- \cdot \overline{\tau}(F)(a)$ is fully faithful because $\overline{\tau}(F)(a)$ is Cartesian, so in particular $\overline{\kappa}(f)$ is an isomorphism if $\overline{\kappa}(f) \cdot \overline{\tau}(F)(a)$ is, and this 2-cell is an isomorphism because every other 2-cell arising in the definition of $\overline{\kappa}(f)$ is an isomorphism. ∎

## A limit lifting theorem for fibrations between bicategories

Let $p : \mathbb{E} \to \mathbb{B}$ be a strict functor between bicategories. For the sake of this paper, the "fiber" of $p$ over $b$ refers the bicategory whose 0-cells are objects lying strictly over $b$, whose 1-cells $f : e \to e'$ are 1-cells of $\mathbb{E}$ equipped with a distinguished 2-isomorphism $p(f) \cong 1_b$, and whose 2-cells are 2-cells of $\mathbb{E}$ respecting the associated 2-isomorphism.

Let $p : \mathbb{E} \to \mathbb{B}$ be a fibration, $f : b \to b'$ a 1-cell in $\mathbb{B}$, and $\mathbb{A}$ a bicategory. Let $F : \mathbb{A} \to p^{-1}(b')$ be a pseudofunctor valued in the fiber over $b'$; we can regard $F$ as a pair $(F_0, \rho)$ where $F_0 : \mathbb{A} \to \mathbb{E}$ and $\rho : F_0 \cdot p \Rightarrow 1_{b'}$ is an invertible icon.

There is a canonical choice of pseudonatural transformation $\tau : b \Rightarrow F \cdot p$ given by the data of the icon $\rho$.

A **reindexing of $F$ along $f$** consists of a functor $\overline{F} : \mathbb{A} \to \mathbb{E}$ lying strictly over the constant functor $1_{b'}$ together with a Cartesian pseudonatural transformation $\overline{\tau} : \overline{F} \Rightarrow F_0$ lying over $\tau$. Evidently, $\overline{F}$ factors through the inclusion of the fiber of $b$ into $\mathbb{E}$.

Suppose $F : A \to p^{-1}(b')$ has a (conical) limit in the fiber. Because of our work to show that $- \cdot p$ is a fibration, it is evident that for any $f : b \to b'$, the limit cone can be reindexed along $f$, and if this reindexed cone is a limit cone, we say that $f$ preserves the limit.

We say that reindexing along $f$ preserves limits of shape $\mathbb{A}$ if it preserves limits for all pseudofunctors $\mathbb{A} \to p^{-1}(b)$.

> **Theorem 0.19** (Limit lifting theorem): Let $p : \mathbb{E} \to \mathbb{B}$ be a fibration of bicategories, and let $\mathbb{A}$ be a diagram bicategory. If $\mathbb{B}$ has pseudo (respectively, oplax) limits of shape $\mathbb{A}$, and the fibers of $p$ have pseudo (respectively, oplax) limits of shape $\mathbb{A}$ preserved by reindexing, then $\mathbb{E}$ has pseudo (respectively, oplax) limits of shape $\mathbb{A}$.

*Proof*: The proof is almost entirely formal, relying on the fibration structure for pseudofunctors and pseudo (respectively, oplax) natural transformations we have established; the same proof works in both the pseudo and oplax case. Let $J : \mathbb{A} \to \mathbb{E}$ be a diagram. Let $(s, \tau : s \Rightarrow J \cdot p)$ be the limit of $J \cdot p$, where $s$ is an object of $\mathbb{B}$ and by abuse of notation also the constant functor at $s$. Let $(\overline{s}, \overline{\tau})$ be a Cartesian lift of $\tau$ with codomain $J$; $\overline{s}$ is a diagram in $\mathbb{E}$ that clearly factors through the fiber inclusion of the fiber $p^{-1}(s) \to \mathbb{E}$. Let $(s^*, \tau^*)$ be the limit of this diagram in $p^{-1}(s)$; we claim $(s^*, \tau^* \cdot \overline{\tau})$ is also the limit of the diagram in $\mathbb{E}$. To see this, let $x$ be an arbitrary object of $\mathbb{E}$ and let $\kappa : x \Rightarrow J$ be a cone. Then $\kappa \cdot p$ is a cone for $J \cdot p$, and so $\kappa$ factors through $\tau$ by some 1-cell $g : p(x) \to s$ up to an invertible modification $\alpha$. Because $\overline{\tau}$ is Cartesian, $\kappa$ factors through $\overline{\tau}$ by a morphism $h : x \Rightarrow \overline{s}$, up to an invertible modification $\overline{\beta}$.

Let $\overline{g}_{\lim}$ be a Cartesian 1-cell $g^* s^* \to s^*$. Now, $h : x \to \overline{s}$ factors through the reindexing $\overline{g}_{\overline{s}} : g^*(\overline{s}) \Rightarrow \overline{s}$ by a vertical cone $h_v : x \Rightarrow g^*(\overline{s})$. Because reindexing preserves limits, there is a unique-up-to-isomorphism vertical 1-cell $h' : x \to g^* s^*$ such that $h' \cdot g^*(\tau^*)$ is isomorphic to $h_v$. Furthermore, $h' \cdot \overline{g}_{\lim}$ is the unique 1-cell (up to isomorphism) such that $(h' \cdot \overline{g}_{\lim}) \cdot \tau^* \cong h$, and such that $(h' \cdot \overline{g}_{\lim}) \cdot \tau^* \cdot \overline{\tau} \cong h \cdot \overline{\tau} \cong \kappa$. (Whiskering with a Cartesian 1-cell is fully faithful.)

It remains to prove that for any two cones $\kappa, \kappa' : x \Rightarrow J$ and any modification $\rho : \kappa \Rightarrow \kappa'$, $\rho$ can be factored through $\tau^* \cdot \overline{\tau}$. The proof of this is substantially the same as for 1-cells, and we omit it. ∎